\documentclass[11pt]{article}

\usepackage[margin=1in]{geometry}
\usepackage{amsmath}
\usepackage{amssymb}
\usepackage{graphicx}
\usepackage{hyperref}

\title{A Scale Invariance Property of PCA}
\author{Jonathan Landy}
\date{\today}

\begin{document}
\maketitle

\begin{abstract}
The PCA algorithm is sensitive to changes in measurement scale. Measuring one
variable of a system in inches rather than centimeters, say, alters both its
principal axes and principal eigenvalues. Although this scale dependence is
generally complicated, we show here that it nevertheless obeys a strict
invariance property: under a continuous scale adjustment, the initial state's
$k$-th largest principal component (ordered by eigenvalue) continuously evolves
into the final state's $k$-th largest principal component, for each $k$. In
this sense, we can say that the modes of PCA are ``order-stable'' with respect
to changes in measurement scale --- see Figure \ref{fig:gaussian-stretch} below
for an example. A special case occurs when scaling along directions that are
orthogonal to some modes --- here, apparent eigenvalue crossings can occur.
However, we show that we can interpret these apparent crossings as cases where
the modes instantaneously swap their orientation, in this way maintaining the
required order stability.
\end{abstract}

\begin{figure}[bp]
\centering
\includegraphics[width=0.95\textwidth]{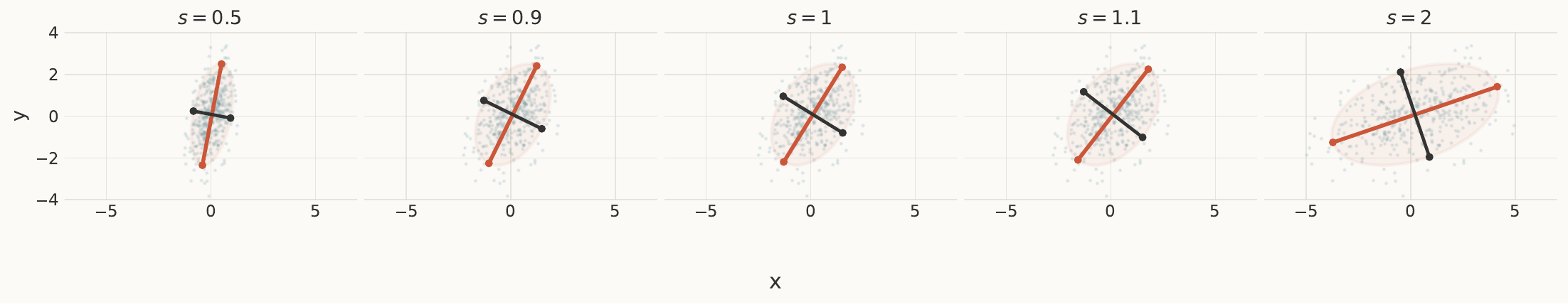}
\caption{PCA response to scaling in 2-d: When we compress along the
$\mathbf{x}$-direction, the smaller (black) principal component rotates to
align with the compression direction. Similarly, when we stretch, the larger
(red) principal component rotates into the scaling direction. At all scales, a
consistent eigenvalue ordering is maintained --- the black mode is always the
smaller of the two and the red is always the larger.}
\label{fig:gaussian-stretch}
\end{figure}

\section{Introduction}

In Principal Component Analysis (``PCA'' for short), we are concerned with the
eigensystem of a given data set's covariance matrix,
\begin{equation}
M = \begin{pmatrix} M_{11} & M_{12} & \ldots \\ M_{21} & M_{22} & \ldots \\ \ldots & & \end{pmatrix}
\label{eq:M}
\end{equation}
The eigensystem provides us with an ellipsoidal approximation to the data set,
giving us a sense of its geometry. Projecting into the space spanned by a
truncated subset of these eigenvectors can often give us a high-fidelity
compression of the data set. These applications make PCA a workhorse algorithm
in data analysis \cite{islp}. The fact that its outputs are sensitive to
measurement scale \cite{miesch1980,wold1987} is therefore a cause of frequent
discomfort: It's almost never clear which choice of units (and which downstream
PCA output) is the objectively ``right'' one for a given application. This
motivates our work here, where we aim to better understand the response of PCA
to changes in measurement scale.

Previous work addressing this issue has considered whether PCA should be
performed on the covariance or correlation matrix, and more generally how
variables should be standardized or scaled before analysis
\cite{miesch1980,wold1987}. Related work has studied the effect of removing
variables, asking when a subset of variables preserves the principal-component
structure of the full data set \cite{jolliffe1972,jolliffe1973}. These studies
consider specific choices of scale, or the limiting case where variables are
removed. Here, we instead ask how the full PCA eigensystem changes continuously
as measurement scales are varied.

To focus our attention, we'll be mostly interested here in the response to a
simple scaling of the first coordinate, writing $x \to s\, x$. In this case,
\begin{equation}
\langle \delta x, \delta y \rangle \to s\, M_{xy}
\label{eq:cov-xy}
\end{equation}
and
\begin{equation}
\langle \delta x^2 \rangle \to s^2\, M_{xx}
\label{eq:cov-xx}
\end{equation}
This sends the covariance matrix to
\begin{equation}
M(s) = \begin{pmatrix} s^2 M_{11} & s M_{12} & s M_{13} & \ldots \\ s M_{21} & M_{22} & M_{23} & \ldots \\ s M_{31} & M_{32} & M_{33} & \ldots \\ & \ldots & & \end{pmatrix}
\label{eq:Ms}
\end{equation}
How do the eigenvectors and eigenvalues of this matrix change as we adjust the
scale factor $s$?

Our key result is that the eigenvalues of \eqref{eq:Ms} never cross as we vary
$s$, but instead always maintain their order. Because this holds for scaling
along any single direction, the result extends to general scalings by
composition, one scaling direction at a time.  The result can be obtained by
transforming the problem to one involving a rank-one perturbation of the
inverse matrix and applying a strict-interlacing theorem for rank-one
perturbations of real symmetric matrices \cite{dobosevych2021}.  Here, we
instead give a direct proof tailored to the positive-definite real symmetric
setting, which also provides a geometric explanation of the result. We first
analyze the ``generic'' case, in three steps:
\begin{itemize}
\item First-order perturbation theory shows that every eigenvalue is strictly
non-decreasing with $s$.
\item In the limit of a very strong compression ($s \to 0$), one eigenvalue
goes to zero, while the other finite eigenvalues approach those of the
subspace orthogonal to the compression direction (here, $\mathbf{x}$). The
Cauchy interlacing theorem constrains these limiting finite eigenvalues to lie
between successive eigenvalues of the original covariance matrix.
\item Combining the last two points via an inductive ``cell'' argument, we
show that --- provided the scaling direction is not orthogonal to any
principal component (the ``generic case'') --- the eigenvalues must maintain
their order under any change in scale.
\end{itemize}

Note that the first two steps here apply in all cases. Only the third requires
the generic assumption. Further, non-orthogonality will hold with probability
$1$ if we select our scaling directions at random --- this is why we call it
the ``generic'' case.  We can also assume that the eigenvalues are
non-degenerate in the generic case. As discussed in the Appendix, when
eigenvalues are degenerate, we can choose an eigenbasis containing an
eigenvector orthogonal to the scaling direction, reducing the problem to the
non-generic case.

We walk through the arguments for the generic case in the next section and
discuss the results in a summary. We cover the non-generic limit and the
orientation swapping effect in an appendix. Some additional results are worked
out in a companion note \cite{landy2026stretch}.

\section{Eigenvalue order maintenance and interlacing}

\subsection{Coupled perturbation system}

To begin our analysis, we'll consider the perturbation equations. If we let $s
\to s + \delta s$ in \eqref{eq:Ms}, the first-order change in $M$ will be
given by
\begin{equation}
\delta M = \begin{pmatrix} 2 s \delta s\, M_{11} & \delta s\, M_{12} & \delta s\, M_{13} & \ldots \\ \delta s\, M_{21} & 0 & 0 & \ldots \\ \delta s\, M_{31} & 0 & 0 & \ldots \\ & \ldots & & \end{pmatrix}
\label{eq:deltaM}
\end{equation}
Applying standard first-order perturbation theory then gives
\begin{equation}
\frac{d\lambda_i}{ds} = \frac{2 \lambda_i \alpha_i^2}{s}
\label{eq:dlambda}
\end{equation}
and
\begin{equation}
\frac{d\vec{v}_i}{ds} = \frac{\alpha_i}{s} \sum_{j \neq i} \frac{\lambda_i +
\lambda_j}{\lambda_i - \lambda_j}\, \alpha_j \vec{v}_j
\label{eq:dv}
\end{equation}
Here, $\vec{v}_i$ is the eigenvector of \eqref{eq:Ms} corresponding to
$\lambda_i$ and $\alpha_i$ is its first component (each evaluated at $s$),
\begin{equation}
\alpha_i \equiv \widehat{e}_1 \cdot \vec{v}_i .
\label{eq:alpha}
\end{equation}
We see that \eqref{eq:dlambda}--\eqref{eq:alpha} give us a coupled set of
differential equations that we can use to solve for the $\vec{v}_i$ and
$\lambda_i$ as we adjust $s$ continuously.

\medskip
\textbf{Key observations}:
\begin{itemize}
\item From \eqref{eq:dlambda}, we see that each non-zero eigenvalue grows
monotonically with $s$, provided its mode is not orthogonal to the scaling
direction.
\item From \eqref{eq:dv}, we see that if a particular eigenvector is
orthogonal to the scaling direction, it will remain fixed for all $s$. In this
case, \eqref{eq:dlambda} implies the eigenvalue will also remain fixed for all
$s$.
\end{itemize}

\subsection{Eigenvalue interlacing under a strong compression}

The Cauchy interlacing theorem \cite{hornjohnson} states the following: Let
$A \in \mathbb{R}^{n\times n}$ be a symmetric matrix, and $B \in
\mathbb{R}^{(n-1) \times (n-1)}$ a principal submatrix of $A$ (obtained by
deleting both the $i$-th row and column of $A$ for some $i$). If $\lambda_1
\leq \lambda_2 \leq \ldots \leq \lambda_n$ are the ordered eigenvalues of $A$
and $\mu_1 \leq \mu_2 \leq \ldots \leq \mu_{n-1}$ are the ordered eigenvalues
of $B$, they interlace:
\begin{equation}
\lambda_{i} \leq \mu_i \leq \lambda_{i + 1}
\label{eq:interlace}
\end{equation}
That is, $\mu_i$ must be at least as large as $\lambda_i$ and it can't be
larger than $\lambda_{i+1}$.

This theorem has direct relevance to our problem in the limit of a strong
compression along a single direction: Qualitatively, what happens in this
limit is that we effectively project the data set into the space orthogonal to
the compression direction (if this is $\mathbf{x}$, say, each data point's $x$
value is pushed to zero). As a result, one eigenvalue will be zero in this
limit (the width in the compression direction), while the other finite
eigenvalues will be those characterizing the space orthogonal to $\mathbf{x}$.
I.e., the finite eigenvalues will match those of the first principal submatrix
of $M$. The Cauchy interlacing theorem therefore applies in this limit.

\medskip
\textbf{Key observations}:
\begin{itemize}
\item Under a strong compression along a single direction, one eigenvalue goes
to zero and $n-1$ eigenvalues remain finite. These new finite eigenvalues
$\mu_i$ and the original eigenvalues $\lambda_i$ satisfy \eqref{eq:interlace}.
\item We detail the asymptotics of the strong compression limit in a
companion note \cite{landy2026stretch}.
\end{itemize}

\subsection{Eigenvalue interlacing (generic case)}

\begin{figure}[htbp]
\centering
\includegraphics[width=0.7\textwidth]{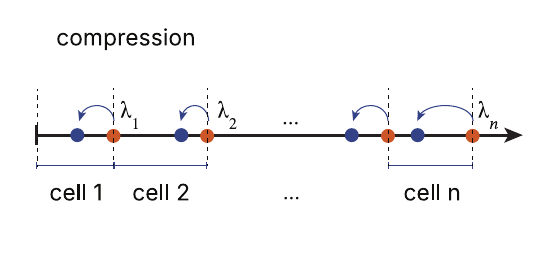}
\caption{Upon compressing by any finite $s < 1$, the eigenvalues must always
monotonically decrease, always maintaining the interlacing property.}
\label{fig:cells}
\end{figure}

In Figure~\ref{fig:cells}, we've illustrated what happens to the eigenvalues
under a generic compression along a single direction. There, we've labeled the
spaces between the original eigenvalues (orange dots) as ordered ``cells''
--- cell $1$ is $(0, \lambda_1)$, and cell $i$ is $(\lambda_{i-1}, \lambda_i)$
for $i \geq 2$. If we compress all the way to $s \rightarrow 0$, our last
section indicates that there must ultimately be exactly one eigenvalue in each
cell (one eigenvalue goes to $0$, while the others interlace with the original
orange eigenvalues).

Consider then what happens under a finite compression. Again, assuming the
generic case, our perturbation equations indicate that each eigenvalue must
strictly and continuously decrease as the compression strength is raised. In
particular, the first eigenvalue $\lambda_1$ must decrease as we start to
compress and it must strictly occupy the first cell throughout (i.e., not
remain at the cell 1 / 2 boundary, its initial position). Under compression, a
second eigenvalue cannot also cross this boundary --- entering the first cell
--- because monotonicity would not allow it to exit before getting to the $s
\to 0$ limit. If this did happen, there would then be two eigenvalues in the
first cell at $s \to 0$, a contradiction. It follows that $\lambda_1$ will
enter cell $1$ under any finite compression, and no second eigenvalue can
enter it as well. An easy inductive argument then implies that $\lambda_2$
must be the only eigenvalue in cell $2$ throughout, and so on.

\begin{figure}[htbp]
\centering
\includegraphics[width=0.7\textwidth]{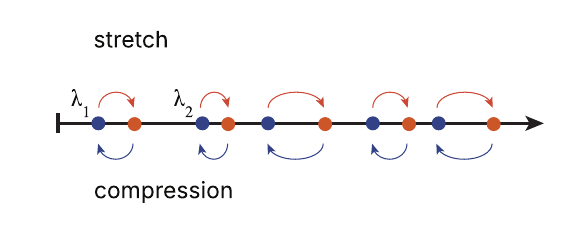}
\caption{Upon a stretch by $s > 1$ along a single direction, the eigenvalues
generically each increase and maintain their ordering. Further, the $i$-th
eigenvalue cannot increase beyond the original $(i+1)$-st eigenvalue.
Similarly, under a generic compression along a single direction, the
eigenvalues each decrease, but without dipping below the original preceding
eigenvalues.}
\label{fig:interlacing}
\end{figure}

To summarize, if we let the eigenvalues of the original (relatively stretched)
matrix be $\lambda_i^{(s)}$ and those of the compressed system be
$\lambda_i^{(c)}$, we have
\begin{equation}
0 < \lambda_1^{(c)} < \lambda_1^{(s)} < \lambda_2^{(c)} < \ldots
< \lambda_n^{(c)} < \lambda_n^{(s)} < \infty,
\label{eq:chain}
\end{equation}
for any generic compression. But of course any compression played in reverse
is a stretch, so \eqref{eq:chain} can also be used to characterize what
happens under a stretch. It follows that upon any finite stretch, the
eigenvalues again interlace, but in the opposite direction: The $i$-th
eigenvalue will increase, but not so much as to cross the original $(i+1)$-st
eigenvalue. Figure~\ref{fig:interlacing} illustrates these results.

\begin{figure}[htbp]
\centering
\includegraphics[width=0.7\textwidth]{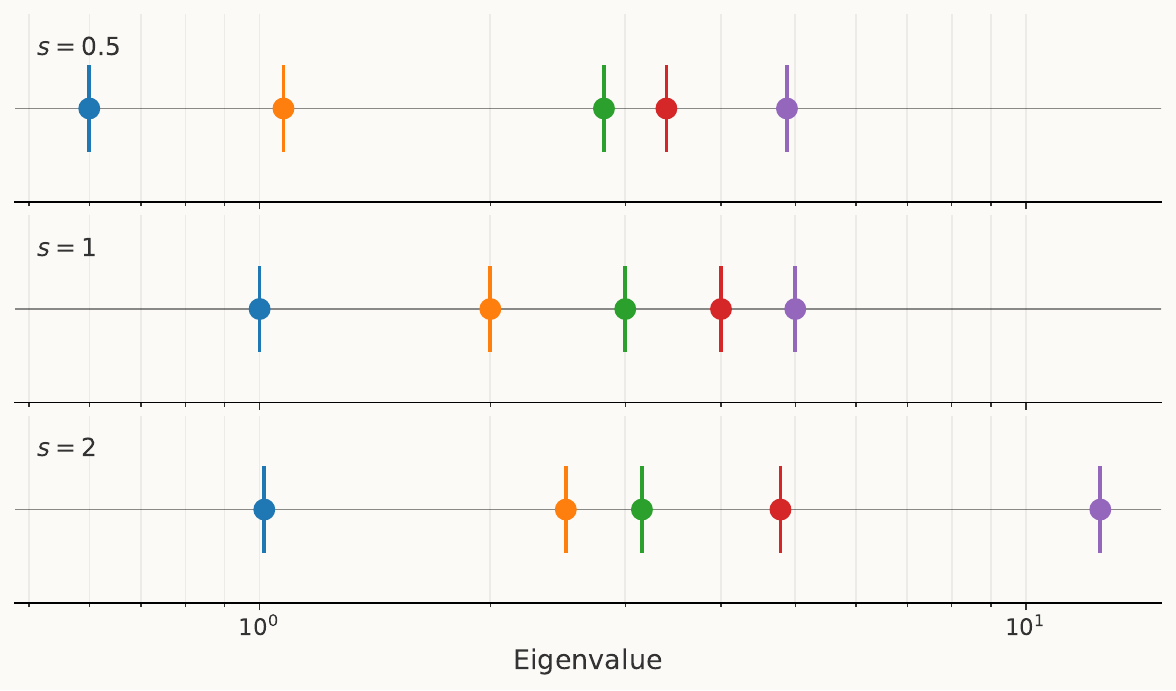}
\caption{Here, we plot the eigenvalues of a random five-dimensional covariance
matrix $M(s)$ at three values of $s$. The eigenvalues never cross one another
and satisfy the strict interlacing property.}
\label{fig:scaled-correlation}
\end{figure}

\medskip
\textbf{Final observations}:
\begin{itemize}
\item Under a generic, single-axis compression or stretch of any finite
strength, the eigenvalues interlace the original eigenvalues and maintain
their initial order. See Figure~\ref{fig:scaled-correlation} for a numerical
example.
\item Main result: Since any set of compressions and stretches can be carried
out via a composition of single-axis compressions and stretches (see our
companion note \cite{landy2026stretch}), it follows that in the generic case
--- where each step's scale direction is not orthogonal to any mode at that
point in the composition --- the modes will also maintain their order. These
results also apply to the non-generic case --- see the Appendix.
\end{itemize}

\section{Summary}

We have shown that monotonicity and asymptotic behavior imply that PCA
eigenvalues cannot cross under continuous changes in scale, and therefore
maintain their ordering. This non-crossing property provides a qualitative
picture of how the principal components evolve under changes in scale. E.g.,
under a strong stretch in one direction, order maintenance tells us that it
must be the largest principal component that rotates into the scaling
direction. Similarly, under a strong compression, it must be the smallest
principal component that rotates into this compression direction --- see
Figure~\ref{fig:gaussian-stretch} for an example. Finally, our cell argument
provides some insight into how the eigenvalues will change under single-axis
scalings, and this can be adapted to understand non-generic cases as well ---
see the Appendix.

Our result may allay some concern during applications that the results of PCA
are sensitive to scale changes. Provided the changes in scale are ``generic''
and not too large, order maintenance implies that the $k$-th mode will
more-or-less maintain its character. However, under non-generic scalings with
apparent crossings, eigenvalue ordering is maintained only through sudden
rotations of the eigenvectors. At such a degeneracy, the two modes effectively
swap their character. This should be kept in mind when considering changes of
scale of this nature (or nearby cases --- again, see the Appendix).

See our companion note \cite{landy2026stretch} for several additional results.
There, we derive the asymptotics of strong single-axis scalings and show that
any two positive-definite covariance matrices can be connected by a finite
sequence of single-axis scalings. That is, the modes of any two
positive-definite covariance matrices are connected in the sense discussed
here. We also analyze the perturbation equations near a ``bounce'' and derive
a simple expression for the change in determinant under scaling.

Although motivated by PCA, these results apply more generally to arbitrary
real symmetric positive-definite matrices, and therefore to a broad class of
spectral problems, including covariance matrices of many mechanical systems.

\appendix
\section{Non-generic scaling and instantaneous rotations}

Here, we'll first illustrate how our cell argument above can be adapted to
anticipate exactly how apparent eigenvalue hopping will occur in the
non-generic regime (again, where the scaling direction is orthogonal to some
modes). Next, we'll consider what happens as we approach this limit from the
generic regime. We'll show that eigenvalue ordering is maintained here via
rapid rotations of the modes. If we interpret the non-generic limit as one
where similar rotations occur instantaneously, we can then formally say that
eigenvalue order maintenance is respected in this regime as well.

Eigenvalue degeneracies also naturally fall within the non-generic regime
considered here. Indeed, at any degeneracy the corresponding eigenspace has
dimension at least two, and therefore necessarily contains a mode orthogonal to
the scaling direction.

\subsection{Non-generic apparent eigenvalue hopping}

\begin{figure}[htbp]
\centering
\includegraphics[width=0.7\textwidth]{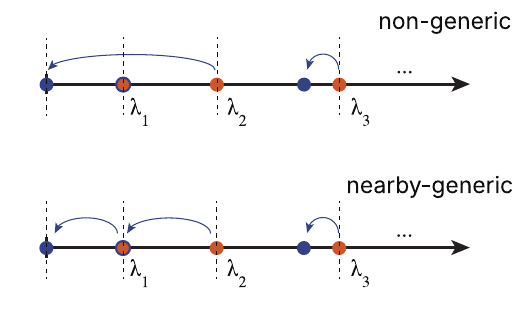}
\caption{In a non-generic case, some eigenvalues will be frozen and can be
hopped. Here, we illustrate what happens under a strong compression orthogonal
to the initial lowest mode. In the nearby-generic case, the eigenvalues do not
cross. Rather, at the near-crossing, $\lambda_2$ bounces into $\lambda_1$,
stopping there. The smaller eigenvalue $\lambda_1$ then continues to zero.}
\label{fig:non-generic}
\end{figure}

The top cell diagram in Figure~\ref{fig:non-generic} provides a simple example
where eigenvalues can be seen to hop in the non-generic limit. There, we
consider a strong compression along a direction orthogonal to the initially
smallest eigenvalue mode (we assume the compression direction is not
orthogonal to any other mode). From our perturbation equations, this means
that $\lambda_1$ must be fixed throughout. However, we also know that in the
limit of a strong compression, one eigenvalue must go to zero and the others
must interlace the original eigenvalues. In this case, since $\lambda_2$ must
decrease, we conclude that $\lambda_1$ and $\lambda_2$ must be the two (and
only two) eigenvalues occupying the first two cells throughout the
compression process. Since one of these must go to zero, it follows that
$\lambda_2 \to 0$ as $s \to 0$, hopping below $\lambda_1$. The logic then
continues as before, with the remaining eigenvalues interlacing --- e.g.,
$\lambda_3$ must be the only eigenvalue in cell $3$, etc.

We note that it's also possible to construct simple examples where hopping of
intermediate eigenvalues occurs, provided these are orthogonal to the scaling
direction. A special case occurs when scaling along a single principal
component. In this case, all eigenvalues are frozen except for that
corresponding to the scaling direction and this can move to any finite value
under an appropriate scaling.

\medskip
\textbf{Key observations}:
\begin{itemize}
\item From the perturbation equations \eqref{eq:dlambda} and \eqref{eq:dv}, we
see that modes orthogonal to scaling directions cannot move. Strong scalings
can then force ``apparent'' crossings.
\item As illustrated here, our cell argument can often be adapted to
understand qualitatively what happens under non-generic scalings.
\end{itemize}

\subsection{Rapid rotations save order maintenance}

The eigensystem of a matrix is continuous in its components, at least away
from any eigenvalue degeneracies. This introduces the question of how we can
break order maintenance under a non-generic scale change, but retain it for a
nearby generic transformation.

The resolution to the issue can be understood through consideration of the
bottom cell diagram of Figure~\ref{fig:non-generic}. Here, we show what
happens for a generic compression very near to the non-generic one considered
above it. In this case, the eigenvalues cannot cross, but must ultimately sit
very close to those of the compressed non-generic system. By continuity, the
eigenvalue locations (ignoring their labels) must approximately match those of
the non-generic system for each point in the compression. What must happen in
this case then is as follows: Before the apparent crossing, $\lambda_1$ will
hardly budge. However, when $\lambda_2$ approaches very close to $\lambda_1$,
$\lambda_1$ will suddenly bounce off and continue towards $0$, leaving
$\lambda_2$ in its prior place. In words, they'll behave like the balls in a
Newton's Cradle toy. This is the only way for the eigenvalues to match in the
two systems throughout the continuous compression process.

\begin{figure}[htbp]
\centering
\includegraphics[width=0.95\textwidth]{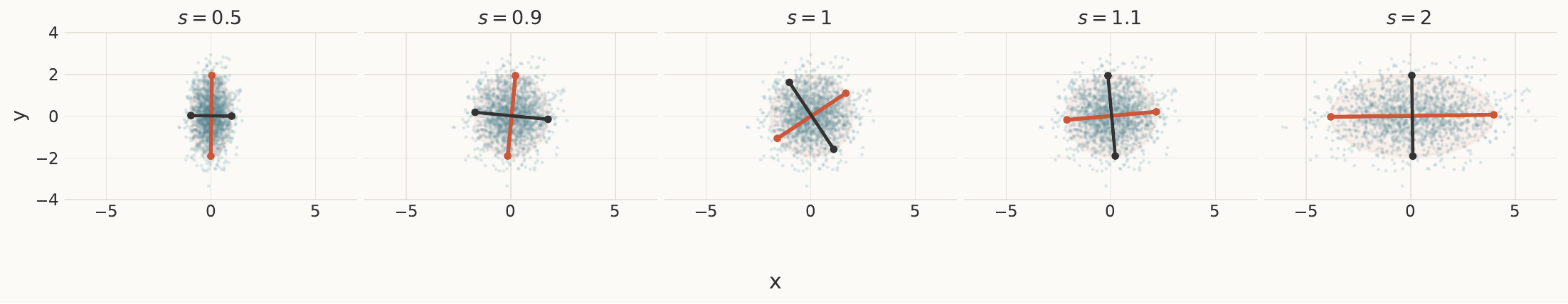}
\caption{A nearly non-generic scaled system in two dimensions. Here, we have
chosen a data set that has one principal component nearly parallel to the
$x$-axis, the direction along which we scale. We show the data and principal
components at various values of $s$, highlighting the rapid rotation that
occurs around $s=1$.  Here, the two eigenvalues approach each other, and we see
a rapid rotation of the modes, preventing an eigenvalue crossing.  This
contrasts with Figure \ref{fig:gaussian-stretch} where we see a more gradual
rotation throughout.}
\label{fig:two-d}
\end{figure}

The eigenvectors must also match throughout. Now, the eigenvector
corresponding to $\lambda_1$ does not change at all in the non-generic
scaling. It must similarly barely change as $\lambda_2$ approaches it in the
nearby-generic transformation. For the two systems to continue to match after
the bounce, the eigenvector of $\lambda_2$ must then align with that of
$\lambda_1$ before the bounce. In other words, the two modes must swap
orientation right at the bounce in the generic system. This is possible as we
have a (near-)degeneracy around the bounce (in a degenerate subspace we are
free to choose any orthogonal basis for the eigensystem --- it is this
flexibility that allows for the rapid rotation). Figure~\ref{fig:two-d}
illustrates this effect in a simple two-dimensional example.

\medskip
\textbf{Key observations}:
\begin{itemize}
\item A non-generic mode is not affected by scaling if it is orthogonal to the
scaling direction. For this reason, it looks like the modes can cross in this
case.
\item However, when considering nearby generic scalings, we see that order
maintenance occurs via rapid orientation swaps of the corresponding modes. We
can interpret the non-generic crossings as limiting forms of this behavior,
also respecting order-stability.
\end{itemize}

For another view of this behavior, see our companion note
\cite{landy2026stretch}, where we work out what happens near a crossing of two
eigenvalues via integration of the perturbation equations.

\end{document}